\documentclass[a4paper]{article}

\usepackage{amsmath}
\usepackage{amssymb}
\usepackage[english]{babel}
\usepackage{amsthm}
\usepackage[latin1]{inputenc}
\usepackage{graphicx}
\usepackage{subfigure}
\usepackage{makeidx}
\usepackage{verbatim}
\usepackage{xcolor}
\usepackage{marginnote}
\usepackage{mathrsfs}
\usepackage{hyperref}
\usepackage{tikz}
\usetikzlibrary{calc}

\newtheorem{theorem}{Theorem}

\theoremstyle{definition}

\def\R{\mathbb{R}}

\def\Z{\mathbb{Z}}

\begin{document}

\title{Lazutkin coordinates and the conjugation problem for billiard maps}
\author{Corentin Fierobe}
\date{}
\maketitle

\begin{abstract}
The conjugation problem for billiard maps conjectures that if two strictly convex billiards have conjugated billiard maps, the billiard tables must be homothetic to each other. We show that if two billiard maps are conjugated, the conjugation diffeomorphism is tangent to a Lazutkin change of coordinates. We also recompute the coefficients in the billiard map expansion as the angle ttends to zero, correcting some incorrect expressions in Lazutkin's computations.
\end{abstract}

\section{Introduction}

Given a strictly convex domain $\Omega\subset\R^2$ with smooth boundary, a subject of interest is the billiard dynamics inside $\Omega$ modeling the evolution of a ray of light reflected at the boundary $\partial\Omega$ according to the law of reflection \textit{angle of incidence equals angle of reflection}.

The billiard dynamics is described by the \textit{billiard map} $T_{\Omega}:X_{\Omega}\to X_{\Omega}$, which is a discrete map acting on the space $X_{\Omega}=\R/|\partial\Omega|\Z\times[0,\pi]$ of pairs $(s,\varphi)$ as follows. Consider an arc-length parametrization $\gamma(s)$ of the boundary $\partial\Omega$. Given a pair $(s_0,\varphi_0)\in X_{\Omega}$, consider the oriented line $\ell$ passing through $\gamma(s_0)$ and making an angle $\varphi_0$ with the tangent line to the boundary at $\gamma(s_0)$. By strict convexity of $\Omega$, $\ell$ intersects $\partial\Omega$ at another point, say $\gamma(s_1)$ with $s_1>s_0$. The line $\ell$ is reflected at that point to a new line $\ell'$, which makes an angle $\varphi_1$ with the tangent line at $\gamma(s_1)$. The billiard map is defined as 
\[
T_{\Omega}(s_0,\varphi_0)=(s_1,\varphi_1).
\]
If the boundary $\partial\Omega$ is $\mathscr C^r$-smooth for a certain $r\geq 2$, the billiard map $T_{\Omega}$ is $\mathscr C^{r-1}$-smooth and is a famous example of \textit{symplectic twist map}, see \cite{FKS_lecturenotes} for a survey on strictly convex billiard maps.
\vspace{0.2cm}

In this paper, we address the famous \textit{conjugation problem} for billiard maps, which can be stated as follows. Let $\Omega_1$ and $\Omega_2$ be two strictly convex billiard domains with respective billiard maps $T_{\Omega_1}$ and $T_{\Omega_2}$.  We say that $T_{\Omega_1}$ and $T_{\Omega_2}$ are \textit{conjugated by a diffeomorphism} $\Phi: X_{\Omega_1}\to X_{\Omega_2}$ if 
\[
T_{\Omega_2} = \Phi\circ T_{\Omega_1} \circ \Phi^{-1}.
\]
If two domains differ by rotations, translations and dilations of the plane, then their corresponding billiard maps are conjugated since their boundaries admits conjugated arc-length parametrizations. The conjugation problem asks \textit{whether the only domains $\Omega_1$ and $\Omega_2$ whose billiard maps are conjugated by a sufficiently smooth diffeomorphism are identical up to rotations, translations and dilations.}

The question can be answered positively if $\Phi = \mathrm{id}$. Indeed, if $\varrho(s)$ is the radius of curvature of the boundary $\partial\Omega$ of a domain $\Omega$ at a point of arc-length coordinate $s$, then the billiard map $T$ in $\Omega$ admits the following expansion as $\varphi\to 0$ \cite[\S14 p. 145]{Lazutkin_book}:
\[
\left\{\begin{aligned}
    s_1 &=& s + 2\varrho(s)\varphi + \mathcal O(\varphi^2)\\
    \varphi_1 & = & \varphi + \mathcal O(\varphi^2)
\end{aligned}\right.
\]
From this it follows that if two billiard maps satisfy $T_{\Omega_1}=T_{\Omega_2}$ then the radii of curvature of the
two domains $\Omega_1$ and $\Omega_2$ coincide and thus the domains are isometric.
For a general diffeomorphism $\Phi$, the answer is less obvious. 
\vspace{0.2cm}

In this paper we show that the conjugation problem is closely related to a change of coordinates introduced by Lazutkin, see \cite[\S14 p. 145]{Lazutkin_book}. Given a domain $\Omega$ with radius of curvature $\varrho$, he introduced a new pair of coordinates $(x,y)$ defined on a subset $Y_{\varrho}\subset\R/\Z\times\R$ by
\[
x = C_{\varrho}\int_0^s\varrho^{-2/3}(\sigma)d\sigma,
\qquad
y = 4C_{\varrho}\varrho^{1/3}(s)\sin\left(\frac{\varphi}{2}\right),
\]
where $C_{\varrho} = \left(\int_0^{|\partial\Omega|}\varrho^{-2/3}(s)ds\right)^{-1}>0$ is a normalization constant. In these coordinates, the billiard map in $\Omega$ admits the following expansion
\[
\left\{\begin{aligned}
    x_1 &=& x+y+\mathcal O(y^3)\\
    y_1 &=& y +\mathcal O(y^4)
\end{aligned}\right.
\]
Observing that the billiard map in these coordinates is close to the integrable map $(x,y)\mapsto (x+y,y)$, he obtained KAM type results \cite{Lazutkin_KAM}. 
\vspace{0.2cm}

Given a domain $\Omega$ of radius of curvature $\varrho$, we introduce the map $L_{\varrho}:X_{\Omega}\to Y_{\varrho}$ corresponding to the Lazutkin change of coordinates and defined by
\[
L_{\varrho}(s,\varphi) = \left(C_{\varrho}\int_0^s\varrho^{-2/3}(\sigma)d\sigma,
4C_{\varrho}\varrho^{1/3}(s) \sin\left(\frac{\varphi}{2}\right)\right).
\]
Note that arc-length parametrizations are not unique. Consequently, the radius of curvature $\varrho$ of the boundary of a domain admits several representations, given by the functions
\[
 s \in \mathbb{R}/|\partial\Omega|\mathbb{Z} \mapsto \varrho(as+b),
\qquad a \in \{-1,1\},\; b \in \mathbb{R}.
\]
These correspond to the different choices of orientation and origin of the arc-length parameter. Each such choice induces a corresponding definition of the Lazutkin change of coordinates.

\begin{theorem}
\label{thm:conjugation_pb}
    Let $\Omega_1$ and $\Omega_2$ be two strictly convex billiard domains with $\mathscr C^r$-smooth boundary and radii of curvature $\varrho_1$ and $\varrho_2$. Assume that their respective billiard maps $T_{\Omega_1}$ and $T_{\Omega_2}$ are conjugated by a $\mathscr C^3$-smooth diffeomorphism $\Phi$. Then, up to a choice of arc-length parametrizations of the boundaries, we have as $\varphi\to 0$,
    \[
    \Phi(s,\varphi) = L_{\varrho_2}^{-1}\circ L_{\varrho_1}(s,\varphi) +\left(\mathcal O(\varphi),\mathcal O(\varphi^2)\right),
    \qquad s\in\R.
    \]
\end{theorem}
\vspace{0.2cm}

Theorem \ref{thm:conjugation_pb} shows that the Lazutkin change of coordinates appears naturally in the conjugation problem, as it determines the conjugating diffeomorphism $\Phi$ to first order. We therefore give precise computations of the expansion of the billiard map as $\varphi\to 0$ (or equivalently $y\to0$). More precisely, if we write the billiard map as
\[
\left\{\begin{array}{ccl}
    s_1 &=& s+\alpha_1(s)\varphi+\ldots\\
    \varphi_1 &=& \varphi+\beta_2(s)\varphi^2+\ldots
\end{array}\right.
\]
in $(s,\varphi)$-coordinates, or as
\[
\left\{\begin{aligned}[c]
    x_1 &=& x+y+\tilde \alpha_3(x)y^3+\ldots\\
    y_1 &=& y +\tilde \beta_4(x)y^4+\ldots 
\end{aligned}\right.
\]
in Lazutkin coordinates, then we compute explicitly the first coefficients. For the usual expansion, they are given by
{\footnotesize
\begin{equation}
    \label{eq:billiard_expansion}
\begin{aligned}[c]
\alpha_1(s)&=2\varrho(s)\\
\alpha_2(s)&=\frac{4}{3}\varrho'(s)\varrho(s)\\
\alpha_3(s)&=\frac{2}{3}\varrho''(s)\varrho(s)^2+\frac{4}{9}\varrho'(s)^2\varrho(s)\\
\alpha_4(s)&=
\frac{4}{15}\varrho^{(3)}(s)\varrho(s)^3
+\frac{28}{45}\varrho(s)^2\varrho'(s)\varrho''(s)\\
&\qquad\quad+\frac{16}{135}\varrho'(s)^3\varrho(s)
-\frac{2}{45}\varrho'(s)\varrho(s)
\end{aligned}
\quad
\begin{aligned}[c]
\beta_2(s)&=-\frac{2}{3}\varrho'(s)\\
\beta_3(s)&=\frac{4}{9}\varrho'(s)^2-\frac{2}{3}\varrho''(s)\varrho(s)\\
\beta_4(s)&=
-\frac{2}{5}\varrho^{(3)}(s)\varrho(s)^2
+\frac{28}{45}\varrho(s)\varrho'(s)\varrho''(s)\\
&\qquad\quad\qquad-\frac{44}{135}\varrho'(s)^3
-\frac{2}{45}\varrho'(s)
\end{aligned}
\end{equation}}
As a remark, we mention that the ones computed in \cite[\S14 p. 145]{Lazutkin_book} are not all correct\footnote{More precisely, the coefficients $\alpha_3$, $\alpha_4$ and $\beta_4$ are incorrect and given by
\[
\begin{aligned}[c]
\alpha_3&=\frac{2}{3}\varrho''\varrho^2+{\bf \frac{4}{3}}\varrho'^2\varrho\\
\alpha_4&=
\frac{4}{15}\varrho^{(3)}\varrho^3
+{\bf\frac{76}{45}}\varrho^2\varrho'\varrho''
+\frac{16}{135}{\varrho'}^3\varrho
-\frac{2}{45}\varrho'\varrho\\
\beta_4&=
-\frac{2}{15}\varrho^{(3)}\varrho^2
{\bf-\frac{44}{45}}\varrho\varrho'\varrho''
-\frac{44}{135}\varrho'^3
-\frac{2}{45}\varrho'.
\end{aligned}
\]
}, which can lead to incorrect results\footnote{See previous versions of this paper.} and motivated some results in this paper. The coefficients of order at most $4$ in $y$ in \eqref{eq:billiard_expansion_Lazutkin} are given by
{\small
\begin{equation}
    \label{eq:billiard_expansion_Lazutkin}
\begin{aligned}[c]
\tilde\alpha_3(x)&=
-\frac{1}{36}\varrho^{-1}\varrho''(x)
+\frac{1}{27}\varrho^{-2}\varrho'(x)^2
+\frac{1}{96C_{\varrho}^2}\varrho^{-2/3}\\
\tilde\alpha_4(x)&=
-\frac{1}{90}\varrho^{-1}\varrho^{(3)}(x)
+\frac{11}{270}\varrho^{-2}\varrho'(x)\varrho''(x)
-\frac{4}{135}\varrho^{-3}\varrho'(x)^3
-\frac{1}{360C_{\varrho}^2}\varrho^{-5/3}\varrho'(x)\\
\tilde\beta_4(x)&=
\frac{1}{180}\varrho^{-1}\varrho^{(3)}(x)
-\frac{11}{540}\varrho^{-2}\varrho'(x)\varrho''(x)
+\frac{2}{135}\varrho^{-3}\varrho'(x)^3
+\frac{1}{720C_{\varrho}^2}\varrho^{-5/3}\varrho'(x).
\end{aligned}
\end{equation}}

The conjugation problem remains widely open. Using a renormalized version of Mather's beta function, \cite{KalKoudj} showed the existence of a countable family of quantities \textit{a la} Marvizi and Melrose which are preserved by conjugacy. Recent developments towards Birkhoff's conjecture implies that the conjugation problem is true in some cases when one of the domains is an ellipse, see \cite{BialyMironov, KS}.

Theorem \ref{thm:conjugation_pb} is proven in Section \ref{sec:proof_main}, and the computations leading to formulas \eqref{eq:billiard_expansion} and \eqref{eq:billiard_expansion_Lazutkin} are detailed in the appendix.

\section{Acknowledgments}

The author is truly grateful to Qiaoling Wei for stimulating discussions, which led to spotting the errors in Lazutkin coefficients \eqref{eq:billiard_expansion}.
He would also like to thank Alfonso Sorrentino for his help and the interesting conversations, but also Klaudiusz Czudek and Jacopo De Simoi for useful discussions. 
He also acknowledges the support of the Italian Ministry of University and Research's PRIN 2022 grant Stability in Hamiltonian dynamics and beyond, as well as the Department of Excellence grant MathMod@TOV (2023-27) awarded to the Department of Mathematics of the University of Rome Tor Vergata.

\section{Proof of theorem \ref{thm:conjugation_pb}}
\label{sec:proof_main}

Let $\Omega_1$ and $\Omega_2$ be two strictly convex billiards with $\mathscr C^r$-smooth boundary and radii of curvature $\varrho_1$ and $\varrho_2$. Assume that their respective billiard maps $T_{\Omega_1}$ and $T_{\Omega_2}$ are conjugated by a $\mathscr C^3$-smooth diffeomorphism $\Phi:X_{\Omega_1}\to X_{\Omega_2}$:
\[
T_{\Omega_2} = \Phi\circ T_{\Omega_1}\circ \Phi^{-1}.
\]
We rewrite previous equality as $T_{\Omega_2}\circ\Phi = \Phi\circ T_{\Omega_1}$ and we expand $T_{\Omega_2}\circ\Phi$ and $\Phi\circ T_{\Omega_1}$ in $\varphi$.

As $\Phi$ is a diffeomorphism, it preserves the boundaries of $X_{\Omega_1}$ and $X_{\Omega_2}$ hence it admits an expansion of the form
\begin{multline*}
\Phi(s,\varphi) = (\Phi_1(s,\varphi),\Phi_2(s,\varphi))\\
= (a_0(s)+a_1(s)\varphi+\mathcal O(\varphi^2),b_1(s)\varphi+b_2(s)\varphi^2+\mathcal O(\varphi^3))
\end{multline*}
where the $a_i$ and $b_j$ are $|\partial\Omega_1|$-periodic and $a_0:\R/|\partial\Omega_1|\Z\to \R/|\partial\Omega_2|\Z$ is a diffeomorphism.
Without loss of generality, up to translating the arc-length parametrization of $\partial\Omega_2$ by a constant, and eventually applying the transformation $s\mapsto-s$, we can assume that $a_0'>0$ and $a_0(0)=0$.

Consider the asymptotic expansions of the billiard maps given by 
\[
T_{\Omega_1}(s,\varphi) = (s+\alpha_1(s)\varphi+\mathcal O(\varphi^2),\varphi+\beta_2(s)\varphi^2+\mathcal O(\varphi^3))
\]
and 
\[
T_{\Omega_2}(s,\varphi) = (s+\overline{\alpha}_1(s)\varphi+\mathcal O(\varphi^2),\varphi+\overline{\beta}_2(s)\varphi^2+\mathcal O(\varphi^3)).
\]
We first expand $T_{\Omega_2}\circ\Phi(s,\varphi) = (s_1,\varphi_1)$ using Taylor expansions. We drop the $s$ in parentheses to ease reading.
\begin{multline}
\label{eq:expansion_billiard_2s}
    s_1 = \Phi_1+\overline{\alpha}_1(\Phi_1)\Phi_2+\mathcal O(\Phi_2^2)
    = a_0+a_1\varphi+b_1(\overline \alpha_1\circ a_0)\varphi+\mathcal O(\varphi^2) =\\ a_0+(a_1+b_1\overline \alpha_1\circ a_0)\varphi+\mathcal O(\varphi^2)
\end{multline}
\begin{multline}
\label{eq:expansion_billiard_2phi}
    \varphi_1 = \Phi_2+\overline{\beta}_2(\Phi_1)\Phi_2^2+\mathcal O(\Phi_2^3)
    = b_1\varphi+b_2\varphi^2+b_1^2(\overline \beta_2\circ a_0)\varphi^2+\mathcal O(\varphi^3) =\\ b_1\varphi+(b_2+b_1^2\overline \beta_2\circ a_0)\varphi^2+\mathcal O(\varphi^3).
\end{multline}
We then expand $\Phi(s,\varphi)$ in the same way.
\begin{multline}
\label{eq:expansion_billiard_1s}
    \Phi_1 = a_0(s_1)+a_1(s_1)\varphi_1+\mathcal O(\varphi^2)
    = a_0+a_0'\alpha_1\varphi+a_1\varphi+\mathcal O(\varphi^2) =\\ a_0+(a_0'\alpha_1+a_1)\varphi+\mathcal O(\varphi^2)
\end{multline}
\begin{multline}
\label{eq:expansion_billiard_1phi}
    \Phi_2 = b_1(s_1)\varphi_1+b_2(s_2)\varphi_1^2+\mathcal O(\varphi_1^3)
    = (b_1+b_1'\alpha_1\varphi)(\varphi+\beta_2\varphi^2)+b_2\varphi^2+\mathcal O(\varphi^3) =\\ b_1\varphi+(b_2+b_1\beta_2+b_1'\alpha_1)\varphi^2+\mathcal O(\varphi^3).
\end{multline}
Identifying the coefficients of $\varphi$ in Equations \eqref{eq:expansion_billiard_2s} and \eqref{eq:expansion_billiard_1s} and the coefficients of $\varphi^2$ in Equations \eqref{eq:expansion_billiard_2phi} and \eqref{eq:expansion_billiard_1phi}, gives the system
\begin{equation}
\label{eq:system_equations_first_order}
\left\{
\begin{array}{rcl}
    a_0'\alpha_1 & = & b_1\overline\alpha_1\circ a_0  \\
    b_1'\alpha_1+b_1\beta_2 & = & b_1^2\overline\beta_2\circ a_0.
\end{array}\right.
\end{equation}
Observe that $b_1$ cannot vanish, as in the first equation $a_0'$ is never $0$ and $\alpha_1$, $\overline\alpha_1$ are strictly positive by nature, and 
\[
b_1 = \frac{a_0'\alpha_1}{\overline\alpha_1\circ a_0}.
\]
Hence dividing second equation by $b_1$, we obtain
\[
\frac{b_1'}{b_1}\alpha_1+\beta_2 = a_0'\alpha_1\frac{\overline\beta_2\circ a_0}{\overline\alpha_1\circ a_0}.
\]
Dividing again by $\alpha_1$ leads to the equation
\[
\frac{b_1'}{b_1} = a_0'\frac{\overline\beta_2\circ a_0}{\overline\alpha_1\circ a_0}-\frac{\beta_2}{\alpha_1}.
\]
Intergating both sides, we obtain the existence of a constant $c\in \R$ such that
\begin{equation}
\label{eq:ln_b1}
\ln(b_1) = \overline A\circ a_0 - A +k
\end{equation}
where $A$ is a primitive of $\beta_2/\alpha_1$ and $\overline A$ a primitive of $\overline \beta_2/\overline\alpha_1$. Using the expressions of $\alpha_1$ and $\beta_2$, we compute that 
\[
\frac{\beta_2}{\alpha_1} = -\frac{\varrho_1'}{3\varrho_1}
\qquad
\frac{\overline\beta_2}{\overline\alpha_1} = -\frac{\varrho_2'}{3\varrho_2}
\]
hence 
\[
A = -\frac{1}{3}\ln\varrho_1
\qquad 
\overline A = -\frac{1}{3}\ln\varrho_2.
\]
Taking the exponential in \eqref{eq:ln_b1} leads to the existence of a constant $K>0$ such that 
\[
 |b_1| = K\left(\frac{\varrho_2\circ a_0}{\varrho_1}\right)^{-\frac{1}{3}}
\]
Denoting by $\varepsilon\in\{\pm 1\}$ and using back the first equation of \eqref{eq:system_equations_first_order}, we obtain
\begin{equation}
\label{eq:der_a0}
a_0' = \frac{\overline\alpha_1\circ a_0}{\alpha_1}b_1 = \varepsilon K\left(\frac{\varrho_2\circ a_0}{\varrho_1}\right)^{\frac{2}{3}}.
\end{equation}
Since we assumed $a_0'>0$, we have $\varepsilon=1$. Now writing 
\[
\ell_{\varrho}(s) = C_{\varrho}\int_0^s\varrho^{-2/3}(\sigma)d\sigma, \qquad s\in\R
\]
Equation \eqref{eq:der_a0} can be reformulated as
\[
(\ell_{\varrho_2}\circ a_0)' = K\ell_{\varrho_1}'.
\]
Intergating, we deduce the existence of a constant $K'\in\R$ such that
\[
C_{\varrho_2}\ell_{\varrho_2}\circ a_0 = KC_{\varrho_1}\ell_{\varrho_1}+K'.
\]
Evaluating it at $0$ gives $K' = C_{\varrho_2}\ell_{\varrho_2}(a_0(0)) = \ell_{\varrho_2}(0) = 0$. Evaluating it at $|\partial\Omega_1|$ gives
\[
K = \frac{C_{\varrho_2}\ell_{\varrho_2}(|\partial\Omega_2|)}{C_{\varrho_1}\ell_{\varrho_1}(|\partial\Omega_1|)} = 1.
\]
As a consequence, 
\[
a_0 = \ell_{\varrho_2}^{-1}\circ\ell_{\varrho_1},
\qquad
b_1 = \frac{4C_{\varrho_1}\varrho_1^{1/3}}{4C_{\varrho_2}\varrho_2^{1/3}\circ a_0}
\]
and the result follows.

%%%%%%%%%%%%%%%%%%%%%%%%%%%%%%%%%%%%%%%%%%%%%%%%%%%%%%%%%%%
\appendix

\section{Expansion of the billiard map}
\label{sec:expansion_billiard}

This section is devoted to the computations of the coefficients given in Equations \eqref{eq:billiard_expansion} and \eqref{eq:billiard_expansion_Lazutkin}. Given a strictly convex domain $\Omega$ with radius of curvature $\varrho$, the billiard map $T_{\Omega}$ satisfies the following \cite[\S14. p. 145]{Lazutkin_book}
\begin{equation}
\label{eq:implicit_billiard_s}
    \int_{s}^{s_1}\sin\left(\varphi-\int_s^{s'}\varrho^{-1}\right)ds' = 0
\end{equation}
and
\begin{equation}
\label{eq:implicit_billiard_phi}
    \varphi+\varphi_1 = \int_s^{s_1}\varrho^{-1}(s')ds'.
\end{equation}
for each pair $(s,\varphi),(s_1,\varphi_1)\in X_{\Omega}$ such that $T_{\Omega}(s,\varphi) = (s_1,\varphi_1)$. Consider two such pairs and write
\[
\delta(\varphi) = s_1-s.
\]
Since the boundary of $\Omega$ is $\mathscr C^?$-smooth, $\delta$ is four times differentiable in $\varphi$. The map 
\[
G(\varphi) = \int_0^1\sin\left(\varphi-\int_s^{t\delta(\varphi)}\varrho^{-1}\right)dt
\]
is then four times differentiable with zero successive derivatives as it follows from \eqref{eq:implicit_billiard_s}. Setting $u(\varphi,t) = \varphi-\int_s^{t\delta(\varphi)}\varrho^{-1}$, we compute successively
\begin{equation}
    \label{eq:der_G}
\begin{aligned}
    G'(0) &= \int_0^1 \partial_{\varphi}u(0,t)dt\\
    G''(0) &= \int_0^1 \partial_{\varphi}^2u(0,t)dt\\
    G^{(3)}(0) &= \int_0^1 \partial_{\varphi}^3u(0,t)-\partial_{\varphi}u(0,t)^3dt\\
    G^{(4)}(0) &= \int_0^1 \partial_{\varphi}^4u(0,t)-6\partial_{\varphi}u(0,t)\partial_{\varphi}^2u(0,t)dt\\
\end{aligned}
\end{equation}
Now considering the map $F(s') = \int_s^{s'}\varrho^{-1}$, we compute its derivatives at $s$ as given by
\begin{equation}
    \label{eq:derivatives_int_rho}
    \begin{aligned}[c]
        F_1 &= F' = \varrho^{-1}\\
        F_2 &= F'' = -\varrho'\varrho^{-2}\\
        F_3 &= F^{(3)} = 2{\varrho'}^2\varrho^{-3}-\varrho''\varrho^{-2}\\
        F_4 &= F^{(4)} = 6\varrho'\varrho''\varrho^{-3}-6{\varrho'}^3\varrho^{-4}-\varrho^{(3)}\varrho^{-2}.
    \end{aligned}
\end{equation}
We can now compute the $\varphi$-derivatives of $u$ as listed below, where the derivatives of $F$ are implicitely supposed to be composed with $t\delta$:
\begin{equation}
\label{eq:der_u}
\begin{aligned}
     \partial_{\varphi}u&=1-t\delta' F_1\\
     \partial_{\varphi}^2u&= -t\delta''F_1-t^2{\delta'}^2F_2\\
     \partial_{\varphi}^3u&= -t\delta^{(3)}F_1-3t^2\delta'\delta''F_2-t^3{\delta'}^3F_3\\
     \partial_{\varphi}^4u&= -t\delta^{(4)}F_1-4t^2\delta'\delta^{(3)}F_2-3t^2{\delta''}^2F_2-6t^3{\delta'}^2\delta''F_3-t^4{\delta'}^4F_4.
\end{aligned}
\end{equation}
Replacing \eqref{eq:der_u} in the vanishing derivatives of $G$ expressed by \eqref{eq:der_G} and integrating in $t$, we obtain the relations
\[
\begin{aligned}
\text{(}\partial_{\phi}u(0,t) = 1-2t\quad\Longleftarrow\text{)}\qquad\qquad 1-\frac{1}{2}\delta'F_1&=0\\
-\frac{1}{2}\delta''F_1-\frac{1}{3}{\delta'}^2F_2&=0\\
-\frac{1}{2}\delta^{(3)}F_1-\delta'\delta''F_2-\frac{1}{4}{\delta'}^3 F_3&=0\\
-\frac{1}{2}\delta^{(4)}F_1-\frac{4}{3}\delta'\delta^{(3)}F_2-{\delta''}^2F_2-\frac{3}{2}{\delta'}^2\delta''F_3-\frac{1}{5}{\delta'}^4F_4+\delta''F_1+\frac{4}{5}{\delta'}^2F_2&=0
\end{aligned}
\]
where all the maps are evaluated at $\varphi=0$. We therefore obtain
\begin{equation}
\label{eq:der_delta}
\begin{aligned}
\delta' &= \frac{2}{F_1}\\
\delta'' &= -\frac{2}{3}\delta'^2\frac{F_2}{F_1}\\
\delta^{(3)} &= -2\delta'\delta''\frac{F_2}{F_1}-\frac{1}{2}{\delta'}^3\frac{F_3}{F_1}\\
\delta^{(4)} &= -\frac{8}{3}\delta'\delta^{(3)}\frac{F_2}{F_1}-2{\delta''}^2\frac{F_2}{F_1}-3{\delta'}^2\delta''\frac{F_3}{F_1}
-\frac{2}{5}{\delta'}^4\frac{F_4}{F_1}+2\delta''+\frac{8}{5}{\delta'}^2\frac{F_2}{F_1}.
\end{aligned}
\end{equation}
Replacing in \eqref{eq:der_delta} each $F_k$ by its expression given in \eqref{eq:derivatives_int_rho}, we deduce the formula of each $\alpha_j$ according to the relations 
\[
\alpha_j = \frac{\delta^{(j)}(0)}{j!},\qquad j=1,2,3,4.
\]
\vspace{0.2cm}

To obtain the expressions of the $\beta_j$, we use Equation \eqref{eq:implicit_billiard_phi}. A Taylor expansion at $s$ gives first
\[
\begin{aligned}
    \varphi+\varphi_1 &= F_1\delta+\frac{F_2}{2}\delta^2+\frac{F_3}{3!}\delta^3+\frac{F_4}{4!}\delta^4+\mathcal O(\delta^5)\\
    &= F_1(\alpha_1\varphi+\alpha_2\varphi^2+\alpha_3\varphi^3+\alpha_3\varphi^4)
    +\frac{F_2}{2}(\alpha_1^2\varphi^2+2\alpha_1\alpha_2\varphi^3+(\alpha_2^2+2\alpha_1\alpha_3)\varphi^4)\\
    &\quad+\frac{F_3}{3!}(\alpha_1^3\varphi^3+3\alpha_1^2\alpha_2\varphi^4)+\frac{F_4}{4!}\alpha_1^4\varphi^4+\mathcal O(\varphi^5).
\end{aligned}
\]
We therefore deduce the following expressions for the coefficients $\beta_j$ in the Taylor expansion of $\varphi_1$ at $0$: 
\[
\begin{aligned}
    \beta_1 &= F_1\alpha_1-1\\
    \beta_2 &= F_1\alpha_2+\frac{1}{2}\alpha_1^2F_2\\
    \beta_3 &= F_1\alpha_3+F_2\alpha_1\alpha_2+\frac{1}{6}F_3\alpha_1^3\\
    \beta_4 &=F_1\alpha_4+\frac{1}{2}F_2(\alpha_2^2+2\alpha_1\alpha_3)+\frac{1}{2}F_3\alpha_1^2\alpha_2+\frac{1}{24}F_4\alpha_1^4.
\end{aligned}
\]
Relacing the $\alpha_j$ by their expressions given in \eqref{eq:billiard_expansion}, and the $F_k$ by their expansions given in\eqref{eq:derivatives_int_rho}, we obtain the expressions of the $\beta_j$ as claimed in \eqref{eq:billiard_expansion}.

\section{Expansion of the billiard map in Lazutkin coordinates}
\label{sec:expansion_billiard_lazutkin}

In this Section we present the omputations which lead to the expressions of $\tilde\alpha_3$, $\tilde\alpha_4$ and $\tilde \beta_4$. To the pairs $(s,\varphi)$ and $(s_1,\varphi_1)$ we associates the corresponding pairs given in Lazutkin coordinates as $(x,y)$ and $(x_1,y_1)$ respectively.
The first step consists in expanding $x_1-x$ in $\varphi$. To simplify we normalize by dividing by $C_{\varrho}$, which we rename $C$ for simplicity:
\begin{multline*}
\frac{x_1-x}{C} = F_1(\alpha_1\varphi+\alpha_2\varphi^2+\alpha_3\varphi^3+\alpha_4\varphi^4)
+\frac{1}{2}F_2(\alpha_1^2\varphi^2+2\alpha_1\alpha_2\varphi^3+(\alpha_2^2+2\alpha_1\alpha_3)\varphi^4)\\
+\frac{1}{6}F_3(\alpha_1^3\varphi^3+3\alpha_1^2\alpha_2\varphi^4)+\frac{1}{24}F_4\alpha_1^4\varphi^4+\mathcal O(\varphi^5)
\end{multline*}
where each function is evaluated at $s$ and $F_j$ denotes the $j$-th derivative of 
\[
F(s') = \int_s^{s'}\varrho^{-2/3}
\]
at $s$. Hence
\begin{equation}
\label{eq:expansion_x_phi}
\frac{x_1-x}{C} = A_1\varphi+A_2\varphi+A_3\varphi^3+A_4\varphi^4+\mathcal O(\varphi^5)
\end{equation}
where
\[
\begin{aligned}
    A_1 &= F_1\alpha_1\\
    A_2 &= F_1\alpha_2+\frac{F_2}{2}\alpha_1^2\\
    A_3 &= F_1\alpha_3+F_2\alpha_1\alpha_2+\frac{1}{6}F_3\alpha_1^3\\
    A_4 &= F_1\alpha_4+\frac{1}{2}F_2\alpha_2^2+F_2\alpha_1\alpha_3+\frac{1}{2}F_3\alpha_1^2\alpha_2+\frac{1}{24}F_4\alpha_1^4.
\end{aligned}
\]
We compute the $F_j$ as
\[
\begin{aligned}
    F_1 &= \varrho^{-2/3}\\
    F_2 &= -\frac{2}{3}\varrho'\varrho^{-5/3}\\
    F_3 &= \frac{10}{9}{\varrho'}^2\varrho^{-8/3}-\frac{2}{3}\varrho''\varrho^{-5/3}\\
    F_4 &= -\frac{80}{27}{\varrho'}^3\varrho^{-11/3}+\frac{10}{3}\varrho'\varrho''\varrho^{-8/3}-\frac{2}{3}\varrho^{(3)}\varrho^{-5/3}
\end{aligned}
\]
where the derivatives are taken with respect to $s$. Note that we can already compute the $A_j$ in terms of $\varrho$ and its successive derivatives in $s$:
\[
\begin{aligned}
    A_1 &= 2\varrho^{1/3}\\
    A_2 &= 0\\
    A_3 &= \frac{4}{27}{\varrho'}^2\varrho^{1/3}-\frac{2}{9}\varrho''\varrho^{4/3}\\
    A_4 &= -\frac{32}{405}{\varrho'}^3\varrho^{1/3}+\frac{8}{45}\varrho'\varrho''\varrho^{4/3}-\frac{8}{45}\varrho{(3)}\varrho^{7/3}-\frac{2}{45}\varrho'\varrho^{1/3}.
\end{aligned}
\]
Now we transform \eqref{eq:expansion_x_phi} into a Taylor expansion in $y$ using the relation
\[
y = 4C\varrho^{1/3}\sin\left(\frac{\varphi}{2}\right)
\]
which implies
\[
\varphi = 2\arcsin\left(\frac{y}{4C\varrho^{1/3}}\right)
=
\frac{y}{2C\varrho^{1/3}}
+\frac{y^3}{192C^4\varrho}+\mathcal O(y^5).
\]
Replacing $\varphi$ by this expression in \eqref{eq:expansion_x_phi} and multiplying by $C_{\varrho}$ we obtain
\[
    x_1-x = \tilde A_1y+\tilde A_3y^3+\tilde A_4y^4+\mathcal O(y^5)
\]
where
\[
\begin{aligned} 
    \tilde A_1 &= \frac{A_1}{2\varrho^{1/3}}=1\\
    \tilde A_3 &= \frac{1}{C^2}\left(\frac{A_1}{192\varrho}+\frac{A_3}{8\varrho}\right)
    = \frac{1}{C^2}\left(\frac{1}{96}\varrho^{-2/3}+\frac{1}{54}{\varrho'}^2\varrho^{-2/3}-\frac{1}{36}\varrho''\varrho^{1/3}\right)\\
    \tilde A_4 &= \frac{1}{C^3}\left(
    -\frac{2}{405}{\varrho'}^3\varrho^{-1}
    +\frac{1}{90}\varrho'\varrho''
    -\frac{1}{90}\varrho\varrho{(3)}
    -\frac{1}{360}\varrho'\varrho^{-1}
    \right).
\end{aligned}
\]
To express $\tilde\alpha_3$ and $\tilde\alpha_4$, it remains to change the coordinate from $s$ to $x$ in the expressions of $\tilde A_3$ and $\tilde A_4$.
\vspace{0.2cm}

To compute the $\tilde\beta_j$, we start with the following equation, which we expand first in $\varphi$:
\begin{equation}
    \label{eq:starting_tilde_beta}
    y_1= 4C\varrho^{1/3}(s_1)\sin\left(\frac{\varphi_1}{2}\right).
\end{equation}
Considering the successive derivatives $\sigma_j$ of the function $\varrho^{1/3}$, we first obtain
\begin{equation}
    \label{eq:expansion_rho_tier}
    \varrho^{1/3}(s_1) = \varrho^{1/3}+\sigma_1\alpha_1\varphi + (\sigma_1\alpha_2+\frac{1}{2}\sigma_2\alpha_1^2)\varphi^2 +(\sigma_1\alpha_3+\sigma_2\alpha_1\alpha_2+\frac{1}{6}\sigma_3\alpha_1^3)\varphi^3+\mathcal O(\varphi^4).
\end{equation}
We also get the following expansion
\begin{multline}
    \label{eq:expansion_sin_phi}
    \sin\left(\frac{\varphi_1}{2}\right) = \frac{1}{2}\varphi_1-\frac{1}{48}\varphi_1^3+\mathcal O(\varphi_1^5)\\
    = \frac{1}{2}\varphi+\frac{1}{2}\beta_2\varphi^2+\left(\frac{1}{2}\beta_3-\frac{1}{48}\right)\varphi^3+\left(\frac{1}{2}\beta_4-\frac{1}{16}\beta_2\right)\varphi^4+\mathcal O(\varphi^5).
\end{multline}
Multiplying \eqref{eq:expansion_rho_tier} and \eqref{eq:expansion_sin_phi} gives 
\[
\varrho^{1/3}(s_1)\sin\left(\frac{\varphi_1}{2}\right)
= B_1\varphi +B_2\varphi^2+B_3\varphi^3+B_4\varphi^4+\mathcal O(\varphi^5)
\]
where we compute 
\[
\begin{aligned}
    B_1 &= \frac{1}{2}\sigma_0\\
    B_2 &=\frac{1}{2}\beta_2\sigma_0+\frac{1}{2}\sigma_1\alpha_1 \\
    B_3 &= \sigma_0\left(\frac{1}{2}\beta_3-\frac{1}{48}\right)+\frac{1}{2}\beta_2\sigma_1\alpha_1+\frac{1}{2}\sigma_1\alpha_2+\frac{1}{4}\sigma_2\alpha_1^2\\
    B_4 &= \sigma_0\left(\frac{1}{2}\beta_4-\frac{1}{16}\beta_2\right)+\sigma_1\alpha_1\left(\frac{1}{2}\beta_3-\frac{1}{48}\right)+\frac{1}{2}\beta_2(\sigma_1\alpha_2+\frac{1}{2}\sigma_2\alpha_1^2)\\&\qquad\qquad\qquad\qquad\qquad+\frac{1}{2}\sigma_1\alpha_3+\frac{1}{2}\sigma_2\alpha_1\alpha_2+\frac{1}{12}\sigma_3\alpha_1^3.
\end{aligned}
\]
 The $\sigma_j$ can be computed in terms of $\varrho$ ad its $s$-derivatives as follows:
 \[
\begin{aligned}
    \sigma_0 &= \varrho^{1/3}\\
    \sigma_1 &= \frac{1}{3}\varrho'\varrho^{-2/3}\\
    \sigma_2 &= \frac{1}{3}\varrho''\varrho^{-2/3}-\frac{2}{9}{\varrho'}^2\varrho^{-5/3}\\
    \sigma_3 &= \frac{10}{27}{\varrho'}^3\varrho^{-8/3}-\frac{2}{3}\varrho'\varrho''\varrho^{-5/3}+\frac{1}{3}\varrho^{(3)}\varrho^{-2/3}.
\end{aligned}
\]
Combining these expressions for the $\sigma_j$'s with the $\alpha_k$'s given in Equations \eqref{eq:billiard_expansion}, we obtain explicitely
\[
\begin{aligned}
    B_1 &= \frac{1}{2}\varrho^{1/3}\\
    B_2 &= 0\\
    B_3 &= -\frac{1}{48}\varrho^{1/3}\\
    B_4 &= \frac{1}{45}\varrho^{(3)}\varrho^{7/3}-\frac{1}{45}\varrho'\varrho''\varrho^{4/3}+\frac{4}{405}{\varrho'}^3\varrho^{1/3}+\frac{1}{180}\varrho'\varrho^{1/3}.
\end{aligned}
\]We can therefore deduce from the expression of $B_1$, $B_2$ and $B_3$ that
\[
y_1-y = 4CB_4\varphi^4+\mathcal O(\varphi^5) = 4CB_4\left(\frac{y}{2C\varrho^{1/3}}\right)^4+\mathcal O(y^5)
= \tilde B_4 y^4+\mathcal O(y^5)
\]
where 
\[
\tilde B_4 = \frac{\varrho^{-4/3}}{4C^3}B_4 = 
\frac{1}{C^3}\left(
\frac{1}{180}\varrho^{(3)}\varrho-\frac{1}{180}\varrho'\varrho''+\frac{1}{405}{\varrho'}^3\varrho^{-1}+\frac{1}{720}\varrho'\varrho^{-1}
\right).
\]
Expressing $\tilde B_4$ in $x$ coordinates gives the expression of $\tilde\beta_4$ in \eqref{eq:billiard_expansion_Lazutkin}.

\end{document}